\documentclass[11pt,a4paper]{article}
\usepackage{a4,amssymb,latexsym,amsmath,color}
\begin{document}

\newcounter{pkt}
\newenvironment{enum}{\setcounter{pkt}{0}
\begin{list}{\rm\alph{pkt})}{\usecounter{pkt}
\setlength{\topsep}{1ex}\setlength{\labelwidth}{0.5cm}
\setlength{\leftmargin}{1cm}\setlength{\labelsep}{0.25cm}
\setlength{\parsep}{-3pt}}}{\end{list}~\\[-6ex]}

\newcounter{punkt2}
\newenvironment{enumbib}{\setcounter{punkt2}{0}
\begin{small}
\begin{list}{\arabic{punkt2}.}{\usecounter{punkt2}
\setlength{\topsep}{1ex}\setlength{\labelwidth}{0.6cm}
\setlength{\leftmargin}{0.6cm}\setlength{\labelsep}{0.25cm}
\setlength{\parsep}{1pt}}}{\end{list}\end{small}}
    
\newcommand{\bdpm}{\begin{displaymath}}
\newcommand{\edpm}{\end{displaymath}}
\newcommand{\beas}{\begin{eqnarray*}}
\newcommand{\eeas}{\end{eqnarray*}}

\newenvironment{bew}{\vspace*{-0.25cm}\begin{sloppypar}\noindent{\it 
Proof.}}{\hfill\qed\end{sloppypar}\vspace*{0.15cm}}

\newtheorem{theo}{Theorem}          
\newtheorem{cor}[theo]{Corollary}
\newtheorem{prop}[theo]{Proposition}
\newtheorem{exam}[theo]{Example}
\newtheorem{rem}[theo]{Remark}
\newtheorem{rems}[theo]{Remarks}
\newtheorem{lem}[theo]{Lemma}

\newcommand{\brm}{\begin{rm}}
\newcommand{\erm}{\end{rm}}

\newcommand{\qed}{\hfill $\Box$}
\newcommand{\des}{{\sf des}}
\newcommand{\exc}{{\sf exc}}

\begin{center}
{\large\bf A GENERALIZATION OF THE SIMION-SCHMIDT BIJECTION\\ 
FOR RESTRICTED PERMUTATIONS}\\[1cm]
Astrid Reifegerste\\
Institut f\"ur Mathematik, Universit\"at Hannover\\ 
Welfengarten 1, D-30167 Hannover, Germany\\
{\it reifegerste@math.uni-hannover.de}\\[0.5cm]
\renewcommand{\thefootnote}{}
\footnotetext{Date: February 28, 2003}
\end{center}
\vspace*{0.5cm}

\begin{footnotesize}
{\sc Abstract.} 
We consider the two permutation statistics which count the distinct pairs 
obtained from the final two terms of occurrences of patterns 
$\tau_1\cdots\tau_{m-2}m(m-1)$ and $\tau_1\cdots\tau_{m-2}(m-1)m$ in a 
permutation, respectively. By a simple involution in terms of permutation 
diagrams we will prove their equidistribution over the symmetric group.
As special case we derive a one-to-one correspondence between permutations which 
avoid each of the patterns $\tau_1\cdots\tau_{m-2}m(m-1)\in{\cal S}_m$ and 
such ones which avoid each of the patterns $\tau_1\cdots\tau_{m-2}(m-1)m\in{\cal 
S}_m$. For $m=3$ this correspondence coincides with the bijection given by Simion and Schmidt 
in their famous paper on restricted permutations.
\end{footnotesize}
\vspace*{1cm}

\centerline{\large{\bf 1}\hspace*{0.25cm}{\sc Introduction}}
\medskip
In recent time much work has been done investigating permutations with restrictions 
on the patterns they contain. Given a permutation $\pi\in{\cal S}_n$ and a 
permutation $\tau\in{\cal S}_m$, an {\it occurrence of $\tau$ in $\pi$} is an 
integer sequence $1\le i_1<i_2<\ldots<i_m\le n$ such that the letters of the word 
$\pi_{i_1}\pi_{i_2}\cdots\pi_{i_m}$ are in the same relative order as the 
letters of $\tau$. In this context, $\tau$ is called a {\it pattern}. If there is no occurrence we say that $\pi$ {\it avoids} 
$\tau$, or alternatively, $\pi$ is {\it $\tau$-avoiding}. We write ${\cal 
S}_n(\tau)$ to denote the set of $\tau$-avoiding permutations in ${\cal S}_n$, 
and more general, ${\cal S}_n(T)$ for the set of all permutations of length $n$
which avoid each pattern of $T$.\\
A fundamental problem concerning pattern-avoiding permutations is that of 
discovering explicit bijections between the sets ${\cal S}_n(T)$ and ${\cal S}_n(T')$ 
provided they have the same number of elements. The first correspondence of that kind was 
presented by Simion and Schmidt \cite{simion-schmidt}. It proves the well-known fact that $|{\cal S}_n(123)|=|{\cal S}_n(132)|$. We will generalize 
this result.\\[2ex]
In \cite{reifegerste1} and \cite{reifegerste2}, we have used the diagram of a 
permutation to study certain forbidden patterns. For a permutation $\pi\in{\cal S}_n$, we obtain its 
{\it diagram} from the $n\times n$ array representation of $\pi$ by shading, 
for each dot, the cell containing it and the squares that are due south and due 
east of it. Each square left unshaded we call a {\it diagram square}. By the 
construction, the connected components of all diagram squares form Young diagrams. 
For a diagram square, its {\it rank} is defined to be the number of dots northwest of it. 
Clearly, connected diagram squares have the same rank.\\[2ex]
In this paper, we use permutation diagrams again. Section 2 demonstrates 
that the last two terms of occurrences of patterns $\tau\in{\cal S}_m$ with 
$\tau_{m-1}\tau_m=m(m-1)$ and $\tau_{m-1}\tau_m=(m-1)m$, respectively, are closely related to the diagram squares. 
It is an essential property of the permutation statistics which count the number of these 
pairs that they have the same distribution over the symmetric group. To prove 
this fact, a bijection is established in Section 3 which respects these 
statistics. In particular, it will be shown that there are as many permutations 
in ${\cal S}_n$ which avoid each pattern $\tau\in{\cal S}_m$ with $\tau_{m-1}=m$ and 
$\tau_m=m-1$ as permutations which avoid each pattern 
$\tau\in{\cal S}_m$ with $\tau_{m-1}=m-1$ and $\tau_m=m$. For $m=3$ the 
correspondence coincides with Simion-Schmidt's bijection. We conclude with some 
final remarks done in Section 4.
\vspace*{4ex}

\centerline{\large{\bf 2}\hspace*{0.25cm}{\sc Diagrams and occurrences of patterns}}
\medskip
For $m\ge 2$ define the pattern sets
\bdpm
A_m=\{\tau\in{\cal S}_m:\tau_{m-1}=m,\,\tau_m=m-1\}\quad\mbox{and}\quad
B_m=\{\tau\in{\cal S}_m:\tau_{m-1}=m-1,\,\tau_m=m\}.
\edpm
For a permutation $\pi\in{\cal S}_n$, denote by ${\sf a}_m(\pi)$ resp.\! ${\sf 
b}_m(\pi)$ the number of distinct pairs $(i,j)$ with $1\le i<j\le n$ such that $i,j$ are 
the final two terms of an occurrence of a pattern belonging to $A_m$ and $B_m$, 
respectively, in $\pi$. If ${\sf a}_m(\pi)=0$ then $\pi$ avoids each 
pattern of $A_m$; analogously, ${\sf b}_m(\pi)=0$ means the avoidance of every 
pattern of $B_m$. Clearly, ${\sf a}_2(\pi)$ counts the inversions in $\pi$ 
while ${\sf b}_2(\pi)$ says how often $\pi$ contains the pattern $12$. (Apart from $m=2$, 
the numbers ${\sf a}_m(\pi)$ and ${\sf b}_m(\pi)$, respectively, do not coincide 
in general with the total numbers of occurrences of patterns contained in $A_m$ or $B_m$.)\\
For example, the occurrences of $1243$ in $\pi=7\:1\:4\:2\:6\:3\:5\in{\cal S}_7$ 
are the sequences $(2,3,5,7)$, $(2,4,5,6)$, and $(2,4,5,7)$;
$(3,4,5,7)$ is the only occurrence of $2143$. Furthermore, $\pi$ contains once 
the pattern $1234$ and avoids $2134$. Hence ${\sf a}_4(\pi)=2$ and ${\sf 
b}_4(\pi)=1$.\\[2ex]
The number ${\sf a}_m(\pi)$ can immediately be read off from the ranked diagram of 
$\pi$.

\begin{prop} \label{prop1}
Let $\pi\in{\cal S}_n$ be a permutation. Then ${\sf a}_m(\pi)$ equals the 
number of diagram squares of rank at least $m-2$. In particular, $\pi$ avoids  
each pattern of $A_m$ if and only if every diagram square is of rank at most $m-3$. 
\end{prop}

\begin{bew}
It follows from the construction that any diagram square $(i,j)$ of rank at 
least $m-2$ corresponds to an occurrence of a pattern of $A_m$ whose final terms 
are just $i,k$ where $\pi_k=j$.   
\end{bew}
\vspace*{1ex}

By definition, the number ${\sf b}_m(\pi)$ counts the number of non-inversions 
on the letters of $\pi$ which are greater than at least $m-2$ letters to their 
left. (Here a pair $(i,j)$ is called a {\it non-inversion} if $i<j$ and $\pi_i<\pi_j$.) 
It is easy to see that all informations about a permutation are encoded 
in the diagram squares of rank at most $m-3$ and the pairs just counted by ${\sf b}_m(\pi)$.

\begin{prop} \label{prop2}
Any permutation $\pi\in{\cal S}_n$ can completely be recovered from the diagram 
squares having rank at most $m-3$ and the pairs $(i,j)$ arising from the final 
two terms of an occurrence of a pattern belonging to $B_m$. 
\end{prop}

\begin{bew}
By the diagram construction, to know all the diagram squares of rank at most 
$m-3$ means to know the positions of all elements $\pi_i$ for which there are at most $m-3$ 
integers $j<i$ with $\pi_j<\pi_i$. As mentioned before, the 
pairs $(i,j)$ obtained from the end of any occurrence of a pattern of $B_m$ 
are exactly the non-inversions on the set of integers $k$ for which $\pi_k$ 
exceeds at least $m-2$ elements on its left. Thereby 
the positions of these elements are also uniquely determined. 
\end{bew}
\vspace*{1ex}

The proof implies the following procedure for constructing a permutation 
from its diagram squares of rank at most $m-3$ and the final terms of 
occurrences of patterns belonging to $B_m$ in $\pi$.\\
For some integer $m\ge2$, let $D$ be a set of squares and $O$ a 
set of integer pairs that are obtained from a permutation $\pi\in{\cal S}_n$ 
as diagram squares of rank at most $m-3$ and occurrence terms as described, respectively. 
First represent the elements of $D$ as white squares in an 
$n\times n$ array, shaded otherwise. Row by row, put a dot in the leftmost shaded square such 
that there is exactly one dot in each column. Each dot having more than $m-3$ 
dots northwest is being deleted. Then arrange the missing dots such that $O$ is precisely the set 
of non-inversions on these dots, that is, the dot contained in the $i$th row 
lies strictly to the left of the dot contained in the $j$th row if and only if 
$(i,j)\in O$.\\
An efficient way to arrange the dots is the following one. Let 
$r_1<r_2<\ldots<r_s$ be the indices of rows containing no dot, and 
$c_1>c_2>\ldots>c_s$ the indices of columns without a dot. Furthermore, let 
$e_i$ be the number of pairs in $O$ whose first component equals $r_i$. For 
$i=1,\ldots,s$, set $c_i'=c_{e_i+1}$, delete $c_{e_i+1}$ from the sequence $c$, 
and renumber the sequence terms. Put the missing dots in the squares 
$(r_i,c_i')$ where $1\le i\le s$. Note that the second component of the pairs 
contained in $O$ is of no importance for this procedure. It suffices to know 
the first components and the number of their appearance in $O$.

\begin{exam}
\brm
Let $\pi=3\:8\:5\:10\:2\:4\:1\:9\:6\:7\in{\cal S}_{10}$ and $m=5$. The 
left-hand array shows the ranked permutation diagram of $\pi$. All the occurrences of 
patterns of $B_5$ end with $(9,10)$. Thus we obtain: 
\begin{center}                                                        
\definecolor{gray1}{gray}{0.85}
\fboxsep0cm
\fboxrule0cm
\unitlength0.3cm
\begin{picture}(42.5,10)
\put(2,9){\fcolorbox{gray1}{gray1}{\makebox(8,1){}}}
\put(7,8){\fcolorbox{gray1}{gray1}{\makebox(3,1){}}}
\put(4,7){\fcolorbox{gray1}{gray1}{\makebox(6,1){}}}
\put(9,6){\fcolorbox{gray1}{gray1}{\makebox(1,1){}}}
\put(1,5){\fcolorbox{gray1}{gray1}{\makebox(9,1){}}}
\put(3,4){\fcolorbox{gray1}{gray1}{\makebox(7,1){}}}
\put(0,3){\fcolorbox{gray1}{gray1}{\makebox(10,1){}}}
\put(8,2){\fcolorbox{gray1}{gray1}{\makebox(2,1){}}}
\put(5,1){\fcolorbox{gray1}{gray1}{\makebox(5,1){}}}
\put(6,0){\fcolorbox{gray1}{gray1}{\makebox(4,1){}}}
\put(2,0){\fcolorbox{gray1}{gray1}{\makebox(1,10){}}}
\put(7,0){\fcolorbox{gray1}{gray1}{\makebox(1,9){}}}
\put(4,0){\fcolorbox{gray1}{gray1}{\makebox(1,8){}}}
\put(1,0){\fcolorbox{gray1}{gray1}{\makebox(1,6){}}}
\put(3,0){\fcolorbox{gray1}{gray1}{\makebox(1,5){}}}
\put(0,0){\fcolorbox{gray1}{gray1}{\makebox(1,4){}}}
\put(5,0){\fcolorbox{gray1}{gray1}{\makebox(1,2){}}}
\put(19,9){\fcolorbox{gray1}{gray1}{\makebox(8,1){}}}
\put(24,8){\fcolorbox{gray1}{gray1}{\makebox(3,1){}}}
\put(21,7){\fcolorbox{gray1}{gray1}{\makebox(6,1){}}}
\put(24,6){\fcolorbox{gray1}{gray1}{\makebox(3,1){}}}
\put(18,4){\fcolorbox{gray1}{gray1}{\makebox(9,2){}}}
\put(17,0){\fcolorbox{gray1}{gray1}{\makebox(10,4){}}}
\put(19,6){\fcolorbox{gray1}{gray1}{\makebox(1,4){}}}
\put(21,6){\fcolorbox{gray1}{gray1}{\makebox(1,1){}}}
\linethickness{0.3pt}
\multiput(0,0)(0,1){11}{\line(1,0){10}}
\multiput(0,0)(1,0){11}{\line(0,1){10}}
\put(2.5,9.5){\circle*{0.35}}
\put(7.5,8.5){\circle*{0.35}}
\put(4.5,7.5){\circle*{0.35}}
\put(9.5,6.5){\circle*{0.35}}
\put(1.5,5.5){\circle*{0.35}}
\put(3.5,4.5){\circle*{0.35}}
\put(0.5,3.5){\circle*{0.35}}
\put(8.5,2.5){\circle*{0.35}}
\put(5.5,1.5){\circle*{0.35}}
\put(6.5,0.5){\circle*{0.35}}
\multiput(0.5,4.5)(0,1){6}{\makebox(0,0)[cc]{\sf\tiny0}}
\multiput(1.5,6.5)(0,1){4}{\makebox(0,0)[cc]{\sf\tiny0}}
\multiput(3.5,6.5)(0,1){3}{\makebox(0,0)[cc]{\sf\tiny1}}
\multiput(4.5,8.5)(1,0){3}{\makebox(0,0)[cc]{\sf\tiny1}}
\multiput(5.5,6.5)(1,0){2}{\makebox(0,0)[cc]{\sf\tiny2}}
\put(8.5,6.5){\makebox(0,0)[cc]{\sf\tiny3}}
\multiput(5.5,2.5)(1,0){2}{\makebox(0,0)[cc]{\sf\tiny5}}
\multiput(17,0)(0,1){11}{\line(1,0){10}}
\multiput(17,0)(1,0){11}{\line(0,1){10}}
\put(19.5,9.5){\color{red}\circle*{0.35}}
\put(24.5,8.5){\color{red}\circle*{0.35}}
\put(21.5,7.5){\color{red}\circle*{0.35}}
\put(25.5,6.5){\circle*{0.35}}
\put(18.5,5.5){\color{red}\circle*{0.35}}
\put(20.5,4.5){\color{red}\circle*{0.35}}
\put(17.5,3.5){\color{red}\circle*{0.35}}
\put(22.5,2.5){\circle*{0.35}}
\put(23.5,1.5){\circle*{0.35}}
\put(26.5,0.5){\circle*{0.35}}
\multiput(32,0)(0,1){11}{\line(1,0){10}}
\multiput(32,0)(1,0){11}{\line(0,1){10}}
\put(34.5,9.5){\color{red}\circle*{0.35}}
\put(39.5,8.5){\color{red}\circle*{0.35}}
\put(36.5,7.5){\color{red}\circle*{0.35}}
\put(41.5,6.5){\circle*{0.35}}
\put(33.5,5.5){\color{red}\circle*{0.35}}
\put(35.5,4.5){\color{red}\circle*{0.35}}
\put(32.5,3.5){\color{red}\circle*{0.35}}
\put(40.5,2.5){\circle*{0.35}}
\put(37.5,1.5){\circle*{0.35}}
\put(38.5,0.5){\circle*{0.35}}
\put(12,5){\makebox(0,0)[cc]{$:$}}
\put(29.5,5){\makebox(0,0)[cc]{$\longrightarrow$}}
\put(42.5,6.5){\makebox(0,0)[cc]{\sf\tiny4}}
\put(42.5,2.5){\makebox(0,0)[cc]{\sf\tiny8}}
\put(42.5,1.5){\makebox(0,0)[cc]{\sf\tiny9}}
\put(42.5,0.5){\makebox(0,0)[cc]{\sf\tiny10}}
\end{picture}

{\footnotesize{\bf Figure 1}\hspace*{0.25cm}Recovering of a permutation}
\end{center}
Red dots represent the elements of $\pi$ which exceed at most two elements on 
their left. 
Note that $(9,10)$ is the only non-inversion on the set of elements represented 
by black dots in the right-hand array. (The sorting routine yields $c'=(10,9,6,7)$ since 
$e=(0,0,1,0)$.)
\erm
\end{exam}
\vspace*{1ex}

\centerline{\large{\bf 3}\hspace*{0.25cm}{\sc The bijection}}
\medskip

The properties of permutation diagrams given in the previous section are 
essential for the construction of a bijection $\Phi_m$ whose aim it will be to 
prove  

\begin{theo} \label{theo4}
We have $|\{\pi\in{\cal S}_n:{\sf a}_m(\pi)=k\}|=|\{\pi\in{\cal S}_n:{\sf b}_m(\pi)=k\}|$ for all $n$ and $k$.
\end{theo}

Let $\pi\in{\cal S}_n$ be a permutation, and let $D$ be the set of its ranked 
diagram squares. Define $\sigma=\Phi_m(\pi)$ to be the permutation whose 
diagram squares of rank at most $m-3$ 
coincide with the elements of $D$ of rank at most $m-3$, and for which the first 
components of the pairs obtained as final terms of an occurrence of a pattern 
of $B_m$ are just the row indices of the squares in $D$ having rank at 
least $m-2$.\\ 
Before we will analyse this map, we give an example.

\begin{exam} \label{exam5}
\brm
Consider $\pi=3\:8\:5\:10\:2\:4\:1\:9\:6\:7\in{\cal S}_{10}$ again. For $m=5$, 
the map $\Phi_m$ takes $\pi$ to the permutation 
$\sigma=3\:8\:5\:9\:2\:4\:1\:6\:10\:7$:   
\begin{center}                                                        
\definecolor{gray1}{gray}{0.85}
\fboxsep0cm
\fboxrule0cm
\unitlength0.3cm
\begin{picture}(25,10)
\put(2,9){\fcolorbox{gray1}{gray1}{\makebox(8,1){}}}
\put(7,8){\fcolorbox{gray1}{gray1}{\makebox(3,1){}}}
\put(4,7){\fcolorbox{gray1}{gray1}{\makebox(6,1){}}}
\put(9,6){\fcolorbox{gray1}{gray1}{\makebox(1,1){}}}
\put(1,5){\fcolorbox{gray1}{gray1}{\makebox(9,1){}}}
\put(3,4){\fcolorbox{gray1}{gray1}{\makebox(7,1){}}}
\put(0,3){\fcolorbox{gray1}{gray1}{\makebox(10,1){}}}
\put(8,2){\fcolorbox{gray1}{gray1}{\makebox(2,1){}}}
\put(5,1){\fcolorbox{gray1}{gray1}{\makebox(5,1){}}}
\put(6,0){\fcolorbox{gray1}{gray1}{\makebox(4,1){}}}
\put(2,0){\fcolorbox{gray1}{gray1}{\makebox(1,10){}}}
\put(7,0){\fcolorbox{gray1}{gray1}{\makebox(1,9){}}}
\put(4,0){\fcolorbox{gray1}{gray1}{\makebox(1,8){}}}
\put(1,0){\fcolorbox{gray1}{gray1}{\makebox(1,6){}}}
\put(3,0){\fcolorbox{gray1}{gray1}{\makebox(1,5){}}}
\put(0,0){\fcolorbox{gray1}{gray1}{\makebox(1,4){}}}
\put(5,0){\fcolorbox{gray1}{gray1}{\makebox(1,2){}}}
\linethickness{0.3pt}
\multiput(0,0)(0,1){11}{\line(1,0){10}}
\multiput(0,0)(1,0){11}{\line(0,1){10}}
\put(2.5,9.5){\circle*{0.35}}
\put(7.5,8.5){\circle*{0.35}}
\put(4.5,7.5){\circle*{0.35}}
\put(9.5,6.5){\circle*{0.35}}
\put(1.5,5.5){\circle*{0.35}}
\put(3.5,4.5){\circle*{0.35}}
\put(0.5,3.5){\circle*{0.35}}
\put(8.5,2.5){\circle*{0.35}}
\put(5.5,1.5){\circle*{0.35}}
\put(6.5,0.5){\circle*{0.35}}
\multiput(0.5,4.5)(0,1){6}{\makebox(0,0)[cc]{\sf\tiny0}}
\multiput(1.5,6.5)(0,1){4}{\makebox(0,0)[cc]{\sf\tiny0}}
\multiput(3.5,6.5)(0,1){3}{\makebox(0,0)[cc]{\sf\tiny1}}
\multiput(4.5,8.5)(1,0){3}{\makebox(0,0)[cc]{\sf\tiny1}}
\multiput(5.5,6.5)(1,0){2}{\makebox(0,0)[cc]{\sf\tiny2}}
\put(8.5,6.5){\makebox(0,0)[cc]{\sf\tiny3}}
\multiput(5.5,2.5)(1,0){2}{\makebox(0,0)[cc]{\sf\tiny5}}
\multiput(15,0)(0,1){11}{\line(1,0){10}}
\multiput(15,0)(1,0){11}{\line(0,1){10}}
\put(17.5,9.5){\circle*{0.35}}
\put(22.5,8.5){\circle*{0.35}}
\put(19.5,7.5){\circle*{0.35}}
\put(23.5,6.5){\circle*{0.35}}
\put(16.5,5.5){\circle*{0.35}}
\put(18.5,4.5){\circle*{0.35}}
\put(15.5,3.5){\circle*{0.35}}
\put(20.5,2.5){\circle*{0.35}}
\put(24.5,1.5){\circle*{0.35}}
\put(21.5,0.5){\circle*{0.35}}
\put(12.5,5){\makebox(0,0)[cc]{$\longrightarrow$}}
\end{picture}

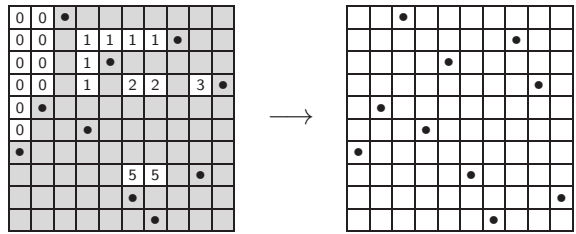
{\footnotesize{\bf Figure 2}\hspace*{0.25cm}Bijection $\Phi_5$, applied to 
$\pi=3\:8\:5\:10\:2\:4\:1\:9\:6\:7$}
\end{center}
The diagram squares having rank at most 2 coincide for $\pi$ and $\sigma$. 
From the row indices of the diagram squares of $\pi$ whose rank is at least 3 
we obtain the first component of each pair arising from the final terms of all 
occurrences of patterns of $B_5$ in $\sigma$. Exactly, we have 
$O=\{(4,*),(8,*),(8,*)\}$. The construction of $\sigma$'s complete array is done 
as described following Proposition \ref{prop2}. With the notations introduced 
there, we have $r=(4,8,9,10)$, $c=(10,9,7,6)$, $e=(1,2,0,0)$, and hence 
$c'=(9,6,10,7)$.
\erm
\end{exam}

As discussed above, the equality of the diagram squares having rank at most $m-3$ 
for $\pi$ and $\sigma=\Phi_m(\pi)$ means that $\sigma_i=\pi_i$ for every $i$ 
for which there are at most $m-3$ integers $j<i$ with $\pi_j<\pi_i$. In 
particular, $\pi$ and $\sigma$ coincide in the $m-2$ first letters. By 
diagram construction, each white square of rank greater than $m-3$ 
is just a pair $(i,\pi_j)$ for which there are at least $m-2$ integers 
$k<i$ with $\pi_k<\pi_j$. Obviously, we have $i<j$ and $\pi_j<\pi_i$. Hence 
both $\pi_i$ and $\pi_j$ are elements exceeding at least $m-2$ elements on 
their left. Consequently, the map $\Phi_m$ is well-defined, and bijective by Proposition 
\ref{prop2} and the remarks done following it.\\[2ex]
It is easy to see that $\Phi_m$ yields the equidistribution of ${\sf a}_m$ and ${\sf 
b}_m$ over the symmetric group. 

\begin{prop}
Let $\pi\in{\cal S}_n$ and $\sigma=\Phi_m(\pi)$, for any $m\ge 2$. Then ${\sf a}_m(\pi)={\sf 
b}_m(\sigma)$.
\end{prop}  

\begin{bew}
By Proposition \ref{prop1}, every pair $(i,j)$ which arises from the two final terms of an 
occurrence of a pattern of $A_m$ in $\pi$ corresponds to a diagram square 
of $\pi$ having rank at least $m-2$, namely $(i,\pi_j)$. It follows immediately from the definition of $\Phi_m$ 
that there is an occurrence of a pattern of $B_m$ in $\sigma$ which ends with 
$(i,k)$ where $k$ depends on $j$.
\end{bew}

\begin{rems}
\brm
\begin{enum}
\item[]
\item By the proof, every occurrence of a pattern of $A_m$ in $\pi$ corresponds in a 
one-to-one fashion to an occurrence of a pattern belonging to $B_m$ in 
$\Phi_m(\pi)$ where both sequences coincide in the $(m-1)$st term. Consequently, 
$\Phi_m$ is even an involution, and we have ${\sf b}_m(\pi)={\sf a}_m(\Phi_m(\pi))$ for all $\pi\in{\cal S}_n$. 
\item The bijection $\Phi_m$ has the advantage of fixing precisely the intersection 
of the sets ${\cal S}_n(A_m)$ and ${\cal S}_n(B_m)$.
\item The map $\Phi_2$ simply takes a permutation $\pi\in{\cal S}_n$ to 
$\sigma\in{\cal S}_n$ with $\sigma_i=n+1-\pi_i$. Note that we have to arrange 
all $n$ dots by the procedure given above. Here $e_i$ equals the number of 
diagram squares in the $i$th row, or, equivalently, the number of integers $j$ 
satisfying $i<j$ and $\pi_i>\pi_j$. The sorting routine yields the permutation 
$\sigma=c'$ whose occurrences of pattern $12$ are just the inversions of $\pi$.
\end{enum}
\erm
\end{rems}

The case $k=0$ of Theorem \ref{theo4} means the {\it Wilf equivalence} of the pattern sets $A_m$ and $B_m$, that is, there are 
as many permutations in ${\cal S}_n$ which avoid each pattern of $A_m$ as those 
which avoid each pattern of $B_m$. An analytical proof of this result was given 
in \cite{mansour}.

\begin{cor}
For each $m\ge 2$, the sets $A_m$ and $B_m$ are Wilf equivalent.
\end{cor}

For $\pi\in{\cal S}_n(A_m)$ the construction of $\sigma=\Phi_m(\pi)$ is even more simple. 
By Proposition \ref{prop1}, every diagram square of $\pi$ is of rank at most $m-3$. Therefore the bijection works as 
follows. Set $\sigma_i=\pi_i$ if there are at most $m-3$ integers $j<i$ satisfying $\pi_j<\pi_i$. 
Then arrange the remaining elements in decreasing order.\\
For example, the permutation 
$\pi=2\:6\:\underline{7}\:1\:\underline{3}\:\underline{4}\:\underline{5}\in{\cal S}_7$ 
avoids both $1243$ and $2143$. We obtain 
$\Phi_4(\pi)=2\:6\:\underline{7}\:1\:\underline{5}\:\underline{4}\:\underline{3}\in{\cal 
S}_7(B_4)$. (All elements which exceed at least two elements on their left are underlined.)\\
Particularly for $m=3$, all the left-to-right mimima of $\pi$ are preserved. (A 
{\it left-to-right minimum} of a permutation $\pi$ is an element $\pi_i$ which 
is smaller than all elements to its left, i.e., $\pi_i<\pi_j$ for every $j<i$.) 
The other elements of $\pi$ are being put on the free positions, in decreasing order. 
This is precisely the description of the bijection between ${\cal S}_n(132)$ 
and ${\cal S}_n(123)$ proposed by Simion and Schmidt in \cite[Prop. 
19]{simion-schmidt}.
\vspace*{4ex}

\centerline{\large{\bf 4}\hspace*{0.25cm}{\sc An open question}}
\medskip
The number of elements of ${\cal S}_n(B_m)$ (or, equivalently, ${\cal S}_n(A_m)$) was determined in \cite{barcucci etal}. 
Similar to the studies of the number of occurrences of certain patterns which are intensified 
in last time, we may ask for the number $|\{\pi\in{\cal S}_n:{\sf 
a}_m(\pi)=k\}|$ for any positive integer $k$. (By Theorem \ref{theo4}, the 
problem of determing $|\{\pi\in{\cal S}_n:{\sf b}_m(\pi)=k\}|$ is equivalent.) As shown in Proposition 
\ref{prop1}, this number just counts the permutations in ${\cal S}_n$ whose diagram 
has exactly $k$ squares of rank at least $m-2$. Here we only consider the special 
case $m=3$ and $k=1$.\\[2ex]
Our proof uses tunnels in Dyck paths which were introduced very recently by 
Elizalde \cite{elizalde}. Recall that a {\it Dyck path} of length $2n$ is a 
lattice path in $\Bbb{Z}^2$ between $(0,0)$ and $(2n,0)$ consisting of up-steps 
$[1,1]$ and down-steps $[1,-1]$ that never falls below the $x$-axis. For any 
Dyck path $d$, a {\it tunnel} is defined to be a horizontal segment between two 
lattice points of $d$ that intersects $d$ only in these two points, and stays 
always below $d$. The {\it length} and {\it height} of a tunnel are measured in 
the lattice. Figure 3 shows a tunnel (drawn with a red line) of length 4 and height 2.
\begin{center}                                                        
\unitlength0.3cm
\begin{picture}(18,5)
\linethickness{0.2pt}
\multiput(0,0)(0,1){6}{\line(1,0){18}}
\multiput(0,0)(1,0){19}{\line(0,1){5}}
\linethickness{0.5pt}
\bezier{150}(0,0)(1.5,1.5)(3,3)
\bezier{150}(3,3)(4.5,1.5)(6,0)
\bezier{50}(6,0)(6.5,0.5)(7,1)
\bezier{50}(7,1)(7.5,0.5)(8,0)
\bezier{100}(8,0)(9,1)(10,2)
\bezier{50}(10,2)(10.5,1.5)(11,1)
\bezier{150}(11,1)(12.5,2.5)(14,4)
\bezier{200}(14,4)(16,2)(18,0)
\linethickness{0.7pt}
\put(12,2){\color{red}\line(1,0){4}}
\end{picture}

{\footnotesize{\bf Figure 3}\hspace*{0.25cm}A tunnel in a Dyck path}
\end{center}

In all Dyck paths of length $2n$ there are ${2n-1\choose n-3}$ tunnels of positive 
height and length at least 4. To verify this, note that there are $nC_n$ 
tunnels in all since every tunnel is associated with an up-step, and the number of Dyck paths of length $2n$ equals the $n$th 
{\it Catalan number} $C_n=\frac{1}{n+1}{2n\choose n}$. The tunnels of height zero correspond 
precisely to {\it returns}, that is, down-steps landing on the $x$-axis. By 
\cite{deutsch}, the total number of returns in Dyck paths of length $2n$ is 
equal to $\frac{3}{2n+1}{2n+1\choose n-1}$. Each tunnel of length 2 and positive 
height is just the connection line of a {\it high peak}. (A {\it high peak} of a Dyck path is an up-step followed 
by a down-step whose common lattice point is at a level greater than $1$.) Their 
number was also given in \cite{deutsch}; it equals ${2n-1\choose n-2}$.

\begin{prop} \label{prop9}
We have $|\{\pi\in{\cal S}_n:{\sf a}_3(\pi)=1\}|={2n-1\choose n-3}$ for all $n$.
\end{prop}

\begin{bew}
We have to count the permutation diagrams having exactly one square of rank 
$r\ge1$. Every such diagram arises from a permutation diagram having only squares of rank zero (note that the 
permutations whose diagram satisfies this condition are just the 
$132$-avoiding ones) by adding the additional white square $(i,j)$ such that $i+j\le 
n+r$, the number of shaded squares $(i,k)$ with $1\le k<j$ is equal to $r$, and the 
number of shaded squares $(k,j)$ with $1\le k<i$ is equal to $r$, too. In 
addition, the $r\times r$ subarray consisting of the squares $(i',j')$ with 
$i-r\le i'\le i-1$ and $j-r\le j'\le j-1$ must be the diagram of any 
permutation in ${\cal S}_r(132)$. This condition is equivalent to 
$i'+j'+r+2\le i+j$ for all diagram squares $(i',j')$ contained in the subarray.\\
Consider now the Dyck path that goes from the lower-left corner to the 
upper-right corner of the array, and travels along the boundary of the connected component 
of all diagram squares of rank zero. 
\begin{center}                               
\definecolor{gray1}{gray}{0.85}
\fboxsep0cm
\fboxrule0cm                         
\unitlength0.3cm
\begin{picture}(7,7)
\put(4,6){\fcolorbox{gray1}{gray1}{\makebox(3,1){}}}
\put(3,5){\fcolorbox{gray1}{gray1}{\makebox(4,1){}}}
\put(3,4){\fcolorbox{gray1}{gray1}{\makebox(2,1){}}}
\put(6,4){\fcolorbox{gray1}{gray1}{\makebox(1,1){}}}
\put(3,3){\fcolorbox{gray1}{gray1}{\makebox(4,1){}}}
\put(2,2){\fcolorbox{gray1}{gray1}{\makebox(5,1){}}}
\put(0,0){\fcolorbox{gray1}{gray1}{\makebox(7,2){}}}
\linethickness{0.3pt}
\multiput(0,0)(0,1){8}{\line(1,0){7}}
\multiput(0,0)(1,0){8}{\line(0,1){7}}
\multiput(0.5,2.5)(0,1){5}{\makebox(0,0)[cc]{\sf\tiny0}}
\multiput(1.5,2.5)(0,1){5}{\makebox(0,0)[cc]{\sf\tiny0}}
\multiput(2.5,3.5)(0,1){4}{\makebox(0,0)[cc]{\sf\tiny0}}
\put(3.5,6.5){\makebox(0,0)[cc]{\sf\tiny0}}
\linethickness{1pt}
\put(0,0){\line(0,1){2}}
\put(0,2){\line(1,0){2}}
\put(2,2){\line(0,1){1}}
\put(2,3){\line(1,0){1}}
\put(3,3){\line(0,1){3}}
\put(3,6){\line(1,0){1}}
\put(4,6){\line(0,1){1}}
\put(4,7){\line(1,0){3}}
\linethickness{0.7pt}
\put(3,4){\color{red}\line(1,0){3}}
\put(6,4){\color{red}\line(0,1){3}}
\put(3,4){\color{red}\line(1,1){3}}
\put(3,4){\circle*{0.3}}
\put(6,7){\circle*{0.3}}
\end{picture}

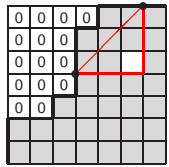
{\footnotesize{\bf Figure 4}\hspace*{0.25cm}Correspondence between diagram 
square of rank $r\ge1$ and Dyck path tunnel}
\end{center}
Any square $(i,j)$ satisfies the above conditions if and only if the line 
that connects the path step contained in the $i$th row (up-step) with the path 
step contained in the $j$th column (down-step) is a tunnel of length 
$2r+2\ge 4$ and height at least 1.
\end{bew}

\begin{rem}
\brm
Thomas \cite{thomas} gives the following alternative combinatorial proof of Proposition 
\ref{prop9} using the permutation statistic ${\bf b}_3$:\\
Let $\pi\in{\cal S}_n$ satisfy ${\sf b}_3(\pi)=1$. Furthermore, let $(i,j)$ be 
the final terms of any occurrence of the pattern $123$ in $\pi$.    Consider now the permutation $\sigma\in{\cal S}_n$ which arises from $\pi$ by exchanging $\pi_i$ 
with $\pi_j$. It is easy to see that $\sigma$ avoids $123$. What can we say 
about the elements $\sigma_i$ and $\sigma_j$? They are successive right-to-left maxima of $\sigma$, 
and there is at least one element to the left of $\sigma_i$ which is smaller 
than $\sigma_j$. (An element is called a {\it right-to-left maximum} of a permutation if it exceeds all the elements on its right.)
In fact, for any $\sigma\in{\cal S}_n(123)$ these two properties characterize 
the pairs $(\sigma_i,\sigma_j)$ whose transposition yields a permutation $\pi$ for which 
${\bf b}_3(\pi)=1$. Consequently, we want to count right-to-left maxima of 
$123$-avoiding permutations for which the set of elements to their right is not 
a complete interval $[1,k]$ for some $k$ or the empty set.\\
In \cite{krattenthaler}, Krattenthaler describes a bijection between $123$-avoiding 
permutations in ${\cal S}_n$ and Dyck paths of length $2n$ having the property that any right-to-left maximum of the kind we 
consider corresponds to a valley at a level greater than zero. (A {\it 
valley} of a Dyck path is a down-step followed by an up-step.) By 
\cite{deutsch}, these valleys are just counted by ${2n-1\choose n-3}$.   
\erm
\end{rem}   

For comparison, Noonan \cite{noonan} proved that the number of permutations 
in ${\cal S}_n$ containing $123$ exactly once is given by $\frac{3}{n}{2n\choose 
n-3}$ while B\'{o}na \cite{bona} showed that there are ${2n-3\choose 
n-3}$ permutations in ${\cal S}_n$ having exactly one $132$-subsequence. By 
\cite[Th. 5.1]{reifegerste1}, the latter permutations are characterized to be 
such ones having exactly one diagram square of rank $1$ and only rank $0$ squares otherwise.
\vspace*{4ex}

\centerline{\large\sc References}
\medskip

\begin{enumbib}

\bibitem{barcucci etal}
E. Barcucci, A. Del Lungo, E. Pergola, and R. Pinzani, 
Permutations avoiding an increasing number of length-increasing forbidden 
subsequences, 
{\it Discrete Math. Theor. Comput. Science} {\bf 4} (2000), 31-44.

\bibitem{bona}
M. B\'{o}na, 
Permutations with one or two $132$-subsequences, 
{\it Discrete Math.} {\bf 181} (1998), 267-274.

\bibitem{deutsch}
E. Deutsch,
Dyck path enumeration,
{\it Discrete Math.} {\bf 204} (1999), 167-202.

\bibitem{elizalde}
S. Elizalde, 
Fixed points and excedances in restricted permutations,
preprint, math.CO/0212221.

\bibitem{krattenthaler}
C. Krattenthaler, 
Permutations with restricted patterns and Dyck paths,
Adv. Appl. Math. 27 (2001) 510-530.

\bibitem{mansour}
T. Mansour, 
Permutations with forbidden patterns, 
Ph.D. Thesis, University of Haifa, 2001.

\bibitem{noonan}
J. Noonan, 
The number of permutations containing exactly one increasing subsequence of 
length three, 
{\it Discrete Math.} {\bf 152} (1996), 307-313.

\bibitem{reifegerste1}
A. Reifegerste, 
On the diagram of $132$-avoiding permutations, 
preprint, math.CO/0208006.

\bibitem{reifegerste2}
A. Reifegerste, 
On the diagram of Schr\"oder permutations, 
{\it Electr. J. Combin.} {\bf 9(2)} (2003), R8.

\bibitem{simion-schmidt}
R. Simion and F. W. Schmidt, 
Restricted Permutations, 
{\it Europ. J. Combin.} {\bf 6} (1985), 383-406.

\bibitem{thomas}
H. Thomas, personal communication.

\end{enumbib}

\end{document}